\documentclass[11pt]{amsart}
\usepackage{amsmath, amsthm, amssymb, amscd, amsfonts}
\usepackage[all]{xy}
\usepackage[OT2,OT1]{fontenc}
\usepackage{newlfont}
\usepackage[dvips]{graphics}
\usepackage[latin1]{inputenc}
\usepackage{eucal}
\usepackage{latexsym}
\usepackage[dvips, dvipsnames, usenames]{color}

\def\wh{\approx_{WH}}

\def\wte{\widetilde{e}}
\def\td{\triangleright}
\def\pop{\cdot_{op}}

\newcommand\cyr{%
\renewcommand\rmdefault{wncyr}%
\renewcommand\sfdefault{wncyss}%
\renewcommand\encodingdefault{OT2}%
\normalfont \selectfont} \DeclareTextFontCommand{\textcyr}{\cyr}

\def\unon{\left\{ 1,\dots ,\theta \right\}}

\newcommand{\map}{{\textcyr{c}}}
\newcommand{\ot}{{\otimes}}
\newcommand{\otv}{{1\le i \le \theta}}
\newcommand{\otvz}{{1\le i, j \le \theta}}
\newcommand{\otrvz}{{1\le i \neq  j \le \theta}}
\newcommand{\qf}{\widetilde{q}}
\newcommand{\qu}{\widehat{q}}
\newcommand{\mf}{m}
\newcommand{\si}{s_{i,F}}
\newcommand{\sE}{s_{i,E}}
\newcommand{\wo}{W_0(\chi)}

\newcommand{\C}{\mathbb C}
\newcommand{\kk}{\Bbbk}
\newcommand{\ku}{\Bbbk}
\newcommand{\N}{\mathbb N}
\newcommand{\Z}{\mathbb Z}

\newcommand{\Sb}{{\mathbb S}}

\def\xx{\mathbb{X}}

\def\bq{\mathfrak{q}}

\def\zt{\Z^{\theta}}

\newcommand{\toba}{{\mathcal B}}
\newcommand{\bB}{\toba}
\newcommand{\G}{{\mathcal G}}

\newcommand{\hlb}{{\mathcal H}}
\newcommand{\Pc}{{\mathcal P}}
\newcommand{\cA}{\mathcal{A}}

\newcommand{\dc}{\Delta(\chi)}

\newcommand\Bg{\mathfrak P}

\newcommand{\e}{\mathbf e}
\newcommand{\f}{\mathbf f}
\newcommand{\ub}{\mathbf u}

\newcommand\ad{\operatorname{ad}_{c}}
\newcommand\gr{\operatorname{gr}}
\newcommand\ord{\operatorname{ord}}
\newcommand\id{\operatorname{id}}
\newcommand\gkd{{\operatorname{GKdim}}}
\newcommand\Aut{\operatorname{Aut}}

\numberwithin{equation}{section} \theoremstyle{plain}
\newtheorem{theorem}{Theorem}[section]
\newtheorem{cor}[theorem]{Corollary}
\newtheorem{lema}[theorem]{Lemma}

\newtheorem{prop}[theorem]{Proposition}

\theoremstyle{definition}
\newtheorem{definition}[theorem]{Definition}

\newtheorem{obs}[theorem]{Remark}

\def\pf{\begin{proof}}
\def\epf{\end{proof}}

\begin{document}

\title[On Nichols algebras with generic braiding]{On Nichols algebras with generic braiding}
\author{Nicol\'as Andruskiewitsch and Iv\'an Ezequiel Angiono}
\address{Facultad of Matem\'atica, Astronom\'\i a y F\'\i sica
\newline \indent
Universidad Nacional of C\'ordoba
\newline
\indent CIEM -- CONICET
\newline
\indent (5000) Ciudad Universitaria, C\'ordoba, Argentina}
\email{andrus@mate.uncor.edu}
 \email{angiono@mate.uncor.edu}
\date{\today}
\thanks{ 2000 {\it Mathematics Subject Classification:} Primary 17B37;
Secondary: 16W20,16W30 \newline \indent {\it Key words and
phrases:} quantized enveloping algebras, Nichols algebras,
automorphisms of non-commutative algebras. }

\begin{abstract} We extend the main result of \cite{AS-crelle} to
braided vector spaces of generic diagonal type using results of
Heckenberger.
\end{abstract}

\maketitle
\section{Introduction}

We fix an algebraically closed field $\kk$ of characteristic 0;
all vector spaces, Hopf algebras and tensor products are
considered over $\kk$. If $H$ is a Hopf algebra, then $G(H):=
\{x\in H -0: \Delta(x) = x\ot x\}$ is a subgroup of the group of
units of $H$; this is a basic invariant of Hopf algebras. Recall
that $H$ is pointed if the coradical of $H$ equals $\kk G(H)$, or
equivalently if any irreducible $H$-comodule is one-dimensional.

The purpose of this paper is to show the validity of the following
classification theorem.

\begin{theorem}\label{fingrowth-lifting}
Let $H$ be a pointed Hopf algebra with finitely generated abelian
group $G(H)$, and generic infinitesimal braiding (see page
\pageref{defi:generic}). Then the following are equivalent:

(a). $H$ is a domain with finite Gelfand-Kirillov dimension.

(b). The group $\Gamma := G(H)$ is free abelian of finite rank,
and there exists a generic datum $\mathcal D$ for $\Gamma$ such
that $H \simeq U({\mathcal D})$ as Hopf algebras.
 \end{theorem}

We refer to the Appendix for the definitions of generic datum and
$U({\mathcal D})$; see \cite{AS-crelle} for a detailed exposition
The general scheme of the proof is exactly the same as for the
proof of \cite[Th. 5.2]{AS-crelle}, an analogous theorem but
assuming ``positive" instead of ``generic" infinitesimal braiding.
The main new feature is the following result.

\begin{lema}\label{rosso}
Let $(V, c)$ be a finite-dimensional braided vector space with
generic braiding. Then the following are equivalent:

(a). $\toba(V)$ has finite Gelfand-Kirillov dimension.

(b). $(V, c)$ is twist-equivalent to a braiding of DJ-type with
finite Cartan matrix.
\end{lema}

Rosso has proved this assuming ``positive" instead of ``generic"
infinitesimal braiding \cite[Th. 21]{R}. Once we establish Lemma
\ref{rosso}, the proofs of Lemma 5.1 and Th. 5.2 in
\cite{AS-crelle} extend immediately to the generic case. Why can
Lemma \ref{rosso} be proved now? Because of the fundamental result
of Heckenberger \cite{H} on Nichols algebras of Cartan type, see
Th. \ref{heckenberger} below. Heckenberger's theorem has as
starting point Kharchenko's theory of PBW-basis in a class of
pointed Hopf algebras. Besides, Heckenberger introduced the
important notion of Weyl groupoid, crucial in the proof of
\cite[Th.  1]{H}.

Here is the plan of this note. In Section \ref{nichols} we
overview Kharchenko's theory of PBW-basis. It has to be mentioned
that related results were announced by Rosso \cite{R2}. See also
\cite{U} for a generalization. We sketch a proof of \cite[L.
19]{R} using PBW-basis. Section \ref{weyl} is devoted to the Weyl
groupoid. We discuss its definition and the proof of \cite[Th.
1]{H} in our own terms. Then we prove Lemma \ref{rosso}.

\section{PBW-basis of Nichols algebras of diagonal type}\label{nichols}

The goal is to describe an appropriate PBW-basis of the Nichols
algebra $\toba(V)$ of a braided vector space of diagonal type. The
argument is as follows. First, there is a basis of the tensor
algebra of a vector space $V$ (with a fixed basis) by Lyndon
words, appropriately chosen monomials on the elements of the
basis. Any Lyndon word has a canonical decomposition as a product
of a pair of smaller words, called the Shirshov decomposition. If
$V$ has a braiding $c$, then for any Lyndon word $l$ there is a
polynomial $[l]_c$ called an hyperletter\footnote{Kharchenko
baptised this elements as superletters, but we suggest to call
them hyperletters to avoid confusions with the theory of
supermathematics.}, defined by induction on the length as a
braided commutator of the hyperletters corresponding to the words
in the Shirshov decomposition. The hyperletters form a PBW-basis
of $T(V)$ and their classes form a PBW-basis of $\toba(V)$.

\subsection{Lyndon words}

\

Let $X$ be a finite set  with a fixed total order: $x_1<\dots
<x_{\theta}$. Let $\xx$ be the corresponding vocabulary-- the set
of words with letters in $X$-- with the lexicographical order.
This order is stable by left, but not by right, multiplication:
$x_1<x_1x_2$ but $x_1x_3>x_1x_2x_3$. However, if $u<v$ and $u$
does not ``begin" by $v$, then $uw<vt$, for all $w,t \in \xx$.

\begin{definition}
A \emph{Lyndon word} is $u \in \xx$, $u\neq 1$, such that $u$ is
smaller than any of its proper ends: if $u=vw$, $v,w \in\xx -
\left\{ 1 \right\}$, then $u<w$. We denote by $L$ the set of
Lyndon words.
\end{definition}

Here are some relevant properties of Lyndon words.

\begin{enumerate}
    \item[(a)] If $u \in L$ and $s \geq 2$, then $u^s \notin L$.
    \item[(b)] Let $u \in \xx-X$. Then, $u$ is Lyndon if and only if for any
    representation $u=u_1 u_2$, with $u_1,u_2 \in \xx$ not empty, one has $u_1u_2=u < u_2u_1$.
    \item[(c)] Any Lyndon word begins by its smallest letter.
    \item[(d)] If $u_1,u_2 \in L, u_1<u_2$, then $u_1u_2 \in L$.
    \item[(e)] If $u,v \in L, u <v$ then $u^k<v$, for all $k \in \N$.
    \item[(f)] Let $u,u_1 \in L$ such that $u=u_3u_2$ and $u_2<u_1$. Then
    $uu_1<u_3u_1$, $uu_1<u_2u_1$.
\end{enumerate}

\begin{theorem}\label{descly} (Lyndon).
Any word $u \in \xx$ can be written in a unique way as a  product
of non increasing Lyndon words: $u=l_1l_2\dots  l_r$, $l_i \in L$,
$l_r \leq \dots \leq l_1$. \qed
\end{theorem}

The \emph{Lyndon decomposition} of $u \in \xx$ is the unique
decomposition given by the Theorem; the $l_i \in L$ appearing in
the decomposition are called the \emph{Lyndon letters} of $u$.

The lexicographical order of $\xx$ turns out to be the same as the
lexicographical order on the Lyndon letters. Namely, if
$v=l_1\dots l_r$ is the Lyndon decomposition of $v$, then $u<v$ if
and only if:
\begin{enumerate}
    \item[(i)] the Lyndon decomposition of $u$ is $u=l_1\dots l_i$, for some $1\leq i <r$, or

    \item[(ii)] the Lyndon decomposition of $u$ is $u=l_1\dots l_{i-1}ll_{i+1}'\dots l_s'$,
    for some $1 \leq i <r$, $s \in \N$ and $l, l_{i+1}',\dots ,l_s'$ in $L$, with $l<l_i$.
\end{enumerate}

Here is another characterization of Lyndon words. See \cite{S, Kh}
for more details.

\begin{theorem}
Let $u \in\xx-X$. Then, $u \in L$ if and only if there exist
$u_1,u_2 \in L$ with $u_1<u_2$ such that $u=u_1u_2$. \qed
\end{theorem}

Let $u \in L-X$. A decomposition $u=u_1u_2$, with $u_1,u_2 \in L$,
is not unique. A very useful decomposition is singled out in the
following way.

\begin{definition} \cite{S}.
Let $u \in L-X$. A decomposition $u=u_1u_2$, with $u_1,u_2 \in L$
such that $u_2$ is the smallest end among those proper non-empty
ends of $u$ is called the \emph{ Shirshov decomposition } of $u$.
\end{definition}

Let $u,v,w \in L$ be such that $u=vw$. Then, $u=vw$ is the
Shirshov decomposition of $u$ if and only if either $v \in X$, or
else if $v=v_1v_2$ is the Shirshov decomposition of $v$, then $w
\leq v_2$.

\subsection{Braided vector spaces of diagonal type}

\

We briefly recall some notions we shall work with; we refer to
\cite{AS-cambr} for more details. A braided vector space is a pair
$(V,c)$, where $V$ is a vector space and $c\in \Aut (V\ot V)$ is a
solution of the braid equation: $(c\otimes \id) (\id\otimes c)
(c\otimes \id) = (\id\otimes c) (c\otimes \id) (\id\otimes c)$. We
extend the braiding to $c:T(V)\ot T(V) \to T(V)\ot T(V)$ in the
usual way. If $x,y\in T(V)$, then the braided bracket is
\begin{equation}\label{eqn:braidedbracket}
[x,y]_c := \text{multiplication } \circ \left( \id - c \right)
\left( x \ot y \right).
\end{equation}

Assume that $\dim V<\infty$ and pick a basis $x_1,\dots
,x_{\theta}$ of $V$, so that we may identify $\kk \xx$ with
$T(V)$. The algebra $T(V)$ has different gradings:

\begin{enumerate}
    \item[(i)] As usual, $T(V) = \oplus_{n\geq 0}T^n(V)$ is $\N_0$-graded.
If $\ell$ denotes the length of a word in $\xx$, then $T^n(V)=
\oplus_{x\in\xx, \, \ell(x) = n}\kk x$.

    \item[(ii)] Let $\e_1, \dots, \e_\theta$ be the canonical basis of $\zt$.
    Then $T(V)$ is also $\zt$-graded, where the degree is
    determined by $\deg x_i = \e_i$, $\otv$.\label{paginatres}
\end{enumerate}

We say that a braided vector space $(V,c)$ is of \emph{diagonal
type} with respect to the basis $x_1, \dots x_\theta$ if there
exist $q_{ij}\in \kk^{\times}$ such that $c(x_i \ot x_j)= q _{ij}
x_j \ot x_i$, $\otvz$. Let $\chi: \zt\times \zt \to \kk^{\times}$
be the bilinear form determined by $\chi(\e_i, \e_j) = q_{ij}$,
$\otvz$. Then
\begin{equation}\label{braiding}
    c(u \ot v)= q _{u,v} v \ot u
\end{equation}
for any $u,v \in \xx$, where $q_{u,v} = \chi(\deg u, \deg v)\in
\kk^{\times}$. Here and elsewhere the degree is with respect to
the $\Z^{\theta}$ grading, see page \pageref{paginatres}. In this
case, the braided bracket satisfies a ``braided" Jacobi identity
as well as braided derivation properties, namely
\begin{align}\label{idjac}
\left[\left[ u, v \right]_c, w \right]_c &= \left[u, \left[ v, w
\right]_c \right]_c
 - q_{u,v} v \ \left[ u, w \right]_c + q_{vw}\left[ u,
w \right]_c \ v,
 \\
\label{der}
    \left[ u,v \ w \right]_c &= \left[ u,v \right]_c w + q_{u,v} v \ \left[ u,w \right]_c,
\\ \label{der2} \left[ u \ v, w \right]_c &= q_{v,w} \left[ u,w \right]_c \ v + u \ \left[ v,w \right]_c,
\end{align}
for any homogeneous $u,v,w \in \toba(V)$.

Let $(V,c)$ be a braided vector space. Let $\widetilde{I}(V)$ be
the largest homogeneous Hopf ideal of the tensor algebra $T(V)$
that has no intersection with $V\oplus {\kk}$. The Nichols algebra
$\toba(V) = T(V)/ \widetilde{I}(V)$ is a braided Hopf algebra with
very rigid properties; it appears naturally in the structure of
pointed Hopf algebras. If $(V,c)$ is of diagonal type, then the
ideal $\widetilde{I}(V)$ is $\zt$-homogeneous hence $\toba(V)$ is
$\zt$-graded. See \cite{AS-cambr} for details. The following
statement, that we include for later reference, is well-known.

\begin{lema}\label{conditions}
\begin{enumerate}
    \item[(a)] If $q_{ii}$ is a root of unit of order $N>1$, then
$x_i^N=0$. In particular, if $\toba(V)$ is an integral domain,
then $q_{ii}=1$ or it is not a root of unit, $i=1,\dots , \theta$.
    \item[(b)] If $i \neq j$, then $(ad_c x_i)^{r}(x_j)=0$ if and only if
    $(r)!_{q_{ii}} \prod _{0 \leq k \leq r-1} (1-q_{ii}^k q_{ij}q_{ji})=0$.
    \item[(c)] If $i\neq j$ and $q_{ij}q_{ji}=q_{ii}^r$, where
    $0\leq -r < \ord(q_{ii})$ (which could be infinite), then $(ad_c
    x_i)^{1-r}(x_j)=0$. \qed
\end{enumerate}
\end{lema}

We shall say that a   braiding $c$ is {\it generic} if it is
diagonal  with matrix $(q_{ij})$ where $q_{ii}$ is not a root of
1, for any $i$.\label{defi:generic}

Finally, we recall that the infinitesimal braiding of a pointed
Hopf algebra $H$ is the braided vector space arising as the space
of coinvariants of the Hopf bimodule $H_1/H_0$, where $H_0 \subset
H_1$ are the first terms of the coradical filtration of $H$. We
refer to \cite{AS-cambr} for a detailed explanation.

\subsection{PBW-basis on the tensor algebra of a braided vector space of diagonal type}

\

We begin by the formal definition of PBW-basis.

\begin{definition}
Let $A$ be an algebra, $P,S \subset A$ and $h: S \mapsto \N \cup
\{ \infty \}$. Let also $<$ be a linear order on $S$. Let us
denote by $B(P,S,<,h)$ the set
\begin{align*}
\big\{ &p\,s_1^{e_1}\dots s_t^{e_t}: t \in \N_0, \quad s_1>\dots
>s_t,\quad s_i \in S, \quad 0<e_i<h(s_i), \quad p \in P \big\}.
\end{align*}

If $B(P,S,<,h)$ is a basis of $A$, then we say that $(P,S,<,h)$ is
a set of \emph{PBW generators} with height $h$, and that
$B(P,S,<,h)$ is a \emph{PBW-basis} of $A$. If $P$, $<$, $h$ are
clear from the context, then we shall simply say that $S$ is a
PBW-basis of $A$.
\end{definition}

\medskip
Let us start with a finite-dimensional braided vector space
$(V,c)$; we fix a basis $x_1,\dots ,x_{\theta}$ of $V$, so that we
may identify $\kk \xx$ with $T(V)$. Recall the braided bracket
\eqref{eqn:braidedbracket}. Let us consider the graded map $\left[
- \right]_c$ of $\kk \xx$ given by
$$
\left[ u \right]_c := \begin{cases} u,& \text{if } u = 1
\text{ or }x \in X;\\
[\left[ v \right]_c, \left[ w \right]_c]_c,  & \text{if } u \in
L, \, \ell(u)>1 \text{ and }u=vw \\ &\qquad\text{ is the Shirshov decomposition};\\
\left[ u_1 \right]_c \dots  \left[ u_t \right]_c,& \text{if } u
\in \xx-L \\ &\qquad
\text{ with Lyndon decomposition  }u=u_1\dots u_t;\\
\end{cases}
$$

Let us now assume that $(V,c)$ is of diagonal type with respect to
the basis $x_1, \dots, x_\theta$, with matrix $(q_{ij})$.

\begin{definition}
Given $l \in L$, the polynomial $\left[ l \right]_c$ is called a
\emph{hyperletter}. Consequently, a \emph{hyperword} is a word in
hyperletters, and a \emph{monotone hyperword} is a hyperword of
the form $W=\left[u_1\right]_c^{k_1}\dots
\left[u_m\right]_c^{k_m}$, where $u_1>\dots >u_m$.
\end{definition}

Let us collect some facts about hyperletters and hyperwords.

\begin{enumerate}
    \item[(a)] Let $u \in L$. Then $\left[ u \right]_c$ is a homogeneous
polynomial with coefficients in $\mathbb{Z} \left[q_{ij}\right]$
and $ \left[ u \right]_c\in u+ \kk \xx^{\ell(u)}_{>u}$.

    \item[(b)]  Given monotone hyperwords  $W,V$, one has
    \[W=\left[w_1\right]_c\dots \left[w_m\right]_c > V=\left[v_1\right]_c\dots  \left[v_t\right]_c, \]
where $w_1 \geq \dots  \geq w_r, v_1 \geq \dots  \geq v_s$, if and
only if     \[w=w_1\dots w_{m} > v=v_i\dots v_t. \] Furthermore,
the principal word  of the polynomial $W$, when  decomposed as sum
of monomials, is $w$ with coefficient 1.
\end{enumerate}

The following statement is due to Rosso.

\begin{theorem} \label{corch} \cite{R2}.
Let $m,n \in L$, with $m<n$. Then $\left[\left[m\right]_c,
\left[n\right]_c \right]_c$ is a  $\mathbb{Z}
\left[q_{ij}\right]$-linear combination of monotone hyperwords
$\left[l_1\right]_c \dots  \left[l_r\right]_c, l_i \in L$, whose
hyperletters satisfy $n>l_i \geq mn$, and such that
$\left[mn\right]_c$ appears in the expansion with non-zero
coefficient. Moreover, for any hyperword
    \[\deg (l_1\dots l_r)= \deg(mn).  \qed\]
 \end{theorem}

The next technical Lemma is crucial in the proof of Theorem
\ref{basehp} below and also in the next subsection. Part (a)
appears in \cite{Kh}, part (b) in \cite{R2}.

\begin{lema}\label{supp}

\begin{enumerate}
    \item[(a)] Any hyperword $W$ is a linear combination  of monotone
    hyperwords bigger than $W$,
    $\left[ l_1 \right] \cdots \left[ l_r \right], l_i \in L$, such that $\deg(W)$ $= \deg(l_1\dots l_r)$,
    and whose hyperletters are between the biggest and the lowest hyperletter of the given word.
    \item[(b)] For any Lyndon  word $l$, let $W_l$ be the vector subspace  of $T(V)$ generated by
    the monotone hyperwords in hyperletters $\left[l_i\right]_c$, $l_i \in L$ such that $l_i \geq l$.
    Then $W_l$ is a subalgebra. \qed
\end{enumerate}
\end{lema}

From this, it can be deduced that the set of monotone hyperwords
is a basis of $T(V)$, or in other words that our first goal is
achieved.

\begin{theorem}\label{basehp} \cite{Kh}.
The set of  hyperletters is a PBW-basis of $T(V)$. \qed
\end{theorem}

\subsection{PBW Basis on quotients of the tensor algebra of a braided vector space of diagonal type}

\

We are next interested in Hopf algebra quotients of $T(V)$. We
begin by describing the comultiplication of hyperwords.

\begin{lema}\label{copro} \cite{R2}.
Let $u \in \xx$, and $u= u_1\dots u_r v^m, \ v, u_i \in L, v<u_r
\leq \dots  \leq u_1$ the Lyndon decomposition of $u$. Then,
\begin{eqnarray*}
        \Delta \left(\left[ u \right]_c\right) &=& 1 \ot \left[ u \right]_c+ \sum ^{m}_{i=0} \binom{ n }{ i } _{q_{v,v}} \left[u_1\right]_c\dots  \left[u_r\right]_c \left[ v \right]_c ^i \ot \left[ v \right]_c^{n-i}
        \\ &+& \sum_{ \substack{ l_1\geq \dots  \geq l_p >l, l_i \in L \\ 0\leq j \leq m } } x_{l_1,\dots ,l_p}^{(j)}
        \ot \left[l_1\right]_c\dots
        \left[l_p\right]_c\left[v\right]_c^j;
\end{eqnarray*}
each $x_{l_1,\dots ,l_p}^{(j)}$ is $\zt$-homogeneous, and
$\deg(x_{l_1,\dots ,l_p}^{(j)})+\deg(l_1\dots  l_p v^j)= \deg(u)$.
\qed
\end{lema}

The following definition appears in \cite{Ha} and is used
implicitly in \cite{Kh}.

\begin{definition}
Let $u,v \in \xx$. We say that $u \succ v$ if and only if either
$\ell(u)<\ell(v)$, or else $\ell(u)=\ell(v)$ and $u>v$
(lexicographical order). This $\succ$ is a total order, called the
\emph{deg-lex order}.
\end{definition}

Note that the empty word 1 is the maximal element for $\succ$.
Also, this order is invariant by right and left multiplication.

\medskip

Let now $I$ be a proper Hopf ideal  of $T(V)$, and set $H=T(V)/I$. Let $\pi: T(V) \rightarrow H$ be the canonical projection. Let us consider the subset of $\xx$:
    \[G_I:= \left\{ u \in \xx: u \notin \\ \xx_{\succ u}+I  \right\}. \]

\begin{prop} \cite{Kh}, see also \cite{R2}.
The set $\pi(G_I)$ is a basis of $H$. \qed
\end{prop}

Notice that

\begin{enumerate}
    \item[(a)] If $u \in G_I$ and $u=vw$, then $v,w \in G_I$.
    \item[(b)] Any word $u \in G_I$  factorizes uniquely as a non-increasing product of Lyndon words in $G_I$.
\end{enumerate}

\noindent Towards finding a PBW -basis the quotient $H$ of $T(V)$,
we look at the set
\begin{equation}\label{setsi}
S_I:= G_I \cap L.
\end{equation}
We then define the function $h_I: S_I \rightarrow \left\{2,3,\dots
\right\}\cup \left\{ \infty \right\}$ by
\begin{equation}\label{defheight}
    h_I(u):= \min \left\{ t \in \N : u^t  \in \kk \xx_{\succ u^t} + I \right\}.
\end{equation}

With these conventions, we are now able to state the main result
of this subsection.

\begin{theorem}\label{basePBW} \cite{Kh}.
Keep the notation above. Then $$B_I':= B\left( \left\{1+I\right\}
, \left[ S_I \right]_c+I, <, h_I \right)$$
    is a PBW-basis of $H=T(V)/I$. \qed
\end{theorem}

The next consequences of the Theorem \ref{basePBW} are used later.
See \cite{Kh} for proofs.

\begin{cor}\label{cor:primero}
A word $u$ belongs to $G_I$ if and only if the corresponding
hyperletter $\left[u\right]_c$  is not a linear combination,
modulo $I$, of greater hyperwords  of the same degree as $u$ and
of hyperwords of lower degree, where all the hyperwords belong to
$B_I$. \qed
\end{cor}

\begin{prop}\label{altf}
In the conditions of the Theorem \ref{basePBW}, if $v \in S_I$ is
such that $h_I(v)< \infty$, then $q_{v,v}$ is a root of unit. In
this case, if $t$ is the order of $q_{v,v}$, then $h_I(v)=t$. \qed
\end{prop}

\begin{cor}\label{cor:segundo}
If $h_I(v):= h < \infty$, then $\left[ v \right]^{h}$ is a linear
combination of monotone hyperwords, in greater hyperletters of
length $h \ell(v)$, and of monotone hyperwords of lower length.
\qed
\end{cor}

\subsection{PBW Basis on the Nichols algebra of a braided vector space of diagonal type}

\

Keep the notation of the preceding subsection. By Theorem
\ref{basePBW}, the Nichols algebra $\toba(V)$ has a PBW-basis
consisting of homogeneous elements (with respect to the
$\zt$-grading). As in \cite{H}, we can even assume that
\begin{itemize}
    \item[$\circledast$] The height of a PBW-generator $\left[ u \right],
\deg(u)=d$, is finite if and only if $2 \leq \ord(q_{u,u}) <
\infty$, and in such case, $h_{\widetilde{I}(V)}(u)=
\ord(q_{u,u})$.
\end{itemize}

This is possible because if the height of $\left[ u \right],
\deg(u)=d$, is finite, then $2 \geq ord(q_{u,u})=m< \infty$, by
Proposition \ref{altf}. And if $2 \leq \ord(q_{u,u})=m < \infty$,
but $h_{\widetilde{I}(V)}(u)$ is infinite, we can add $\left[ u
\right]^m$ to the PBW basis: in this case,
$h_{\widetilde{I}(V)}(u)= \ord(q_{u,u})$, and
$q_{u^m,u^m}=q_{u,u}^{m^2}=1$.

Let $\Delta^+(\toba(V))\subseteq \N^n$ be the set of degrees of
the generators of the PBW-basis, counted with their multiplicities
and let also $\Delta(\toba(V))= \Delta^+(\toba(V)) \cup \left(-
\Delta^+(\toba(V))\right)$. We now show that $\Delta^+(\toba(V))$
is independent of the choice of the PBW-basis with the property
$\circledast$, a fact repeatedly used in \cite{H}.

\medbreak Let $R:= \kk [x_1^{\pm 1}, \ldots, x_{\theta}^{\pm 1}]$,
resp.  $\widehat{R}:= \kk [[x_1^{\pm 1}, \ldots, x_{\theta}^{\pm
1}]]$, the algebra of Laurent polynomials in $\theta$ variables,
resp. formal Laurent series in $\theta$ variables. If $n=(n_1,
\ldots , n_{\theta}) \in \zt$, then we set $X^{n}=X_1^{n_1} \cdots
X_{\theta}^{n_{\theta}}$.  If $T \in \Aut (\zt)$, then we denote
by the same letter $T$ the algebra automorphisms $T: R \rightarrow
R$, $T: \widehat{R} \rightarrow \widehat{R}$, $T(X^n) = X^{T(n)}$,
for all $n \in \zt$. We also set
\begin{eqnarray*}
\bq _h(\mathrm{t}) := \frac{\mathrm{t}^h-1}{\mathrm{t}-1} \in \kk
[\mathrm{t}], \quad h \in \N; \quad \bq _{\infty}(\mathrm{t}):=
\frac{1}{1-\mathrm{t}}= \sum_{s=0}^{\infty} \mathrm{t}^s \in \kk
[[\mathrm{t}]].
\end{eqnarray*}

We say that a $\zt$-graded vector space $V= \oplus _{n \in \zt}
V^n$ is \emph{locally finite} if $\dim V^n < \infty$, for all $n
\in \zt$. In this case, the Hilbert or Poincar\'e series of $V$ is
$\hlb_V = \sum_{n\in \zt} \dim V^n X^n$. If $V$, $W$ are $\zt$
graded, then $ V\ot W = \oplus _{n \in \zt} \left( \oplus_{p \in
\zt} V^p \ot W^{n-p} \right)$ is $\zt$-graded. If $V,W$ are
locally finite and additionally $V^n=W^n=0$, for all $n<M$, for
some $M \in \zt$, then $V\ot W$ is locally finite, and $\hlb_{V
\ot W}=\hlb_V \hlb_W$.

\begin{lema}\label{series}
Let $\chi: \zt \times \zt \rightarrow \kk^{\times}$ be a bilinear
form and set $q_{\alpha}:= \chi(\alpha,\alpha)$, $h_{\alpha}:=
\ord q_{\alpha}$, $\alpha \in \zt$.  Let $N,M \in \N$ and
$\alpha_1,\ldots, \alpha_N, \beta_1,\ldots, \beta_M \in
\N_0^{\theta} \setminus \{ 0 \}$ such that
\begin{equation}\label{eqn:seriesiguales}
 \prod_{1\le i\le N}
\bq_{h_{\alpha_i}}( X^{\alpha_i}) = \prod_{1\le j\le M}
\bq_{h_{\beta_j}}(X^{\beta_j}). \end{equation} Then $N = M$ and
exists $\sigma \in \Sb_N$ such that $\alpha_i =
\beta_{\sigma(i)}$.
\end{lema}

\pf If $\gamma \in \N_0^{\theta} \setminus \{ 0 \}$, then set
\[ C_{\gamma}:= \left\{ (s_1,\ldots,s_N): \sum_{i=1}^N s_i\alpha_i
= \gamma, 0 \leq s_i<h_{\alpha_i}, \quad 1 \leq i \leq N \right\},
\] $c_{\gamma}:= \# C_{\gamma} \in \N_0$. Then the series in
\eqref{eqn:seriesiguales} equals $ 1+\sum_{\N_0^{\theta} \setminus
\{ 0 \}} c_{\gamma} X^{\gamma}$. Let $m_1:= \min \{ |\gamma|:
c_{\gamma} \neq 0 \} $. Then $m_1=c_{\gamma_1}$, for some
$\gamma_1 \in \N_0^{\theta} \setminus \{ 0 \}$, and
$s=(s_1,\ldots,s_N) \in \C_{\gamma_1}$ should belong to the
canonical  basis. Let $I := \{i: \alpha_i= \gamma_1 \}\subseteq \{
1, \ldots N\}$, $J := \{j: \beta_j= \gamma_1 \} \subseteq \{ 1,
\ldots M\}$. Since $c_{\gamma_1}=\# I = \# J$, there exists a
bijection from $I$ to $J$, and moreover, $ \prod_{1\le i\le N, i
\notin I} \bq_{h_{\alpha_i}}( X^{\alpha_i}) = \prod_{1\le j\le M,
j \notin J} \bq_{h_{\beta_j}}(X^{\beta_j})$. The Lemma then
follows by induction on $k=\min\{N,M\}$. \epf

Hence, if $V= \oplus _{n \in \zt} V^n$ is locally finite, and
$\hlb_V = \prod_{1\le i\le N} \bq_{h_{\alpha_i}}( X^{\alpha_i})$,
then the family $\alpha_1,\ldots, \alpha_N$ is unique up to a
permutation.
\medskip

We now sketch a proof of \cite[Lemma 19]{R} using the PBW-basis;
see \cite{Ang} for a complete exposition.

\begin{lema}\label{rossito}
Let $V$ be a braided vector space of diagonal type with matrix
$q_{ij} \in \kk^{\times}$. If $\toba(V)$ is a domain and its
Gelfand-Kirillov dimension is finite, then for any pair $i,j \in
\unon, i \neq j$, there exists $m_{ij}\geq 0$ such that
    \[(ad_cx_i)^{m_{ij} + 1}(x_j)=0. \]
\end{lema}

\pf Let $i \in \unon$. By Lemma \ref{conditions}, either $q_{ii}$
is not a root of the unit or else it is $1$.  Suppose that there
are $i\neq j \in \unon$ such that $(ad_cx_i)^m(x_j) \neq 0$, for
all $m
>0$; say, $i=1$, $j=2$. Hence $q_{11}^mq_{12}q_{21} \neq 1$,
$m \in \N$. Then one can show by induction that $x_1^mx_2 \in
S_I$, using Corollary \ref{cor:primero}. Next, assume that $\left[
x_1^m x_2\right]_c$ has finite height. Then, using Corollary
\ref{cor:segundo}, necessarily $\left[ x_1^mx_2 \right]_c^k=0$ for
some $k$. Let $c_k= \prod^{k-1}_{p=0}(1-q_{11}^pq_{12}q_{21})$.
Using skew-derivations, we obtain that
    \[q_{21}^{-k\frac{(s+1)s}{2}}(s)_{q_{22}^{-1}}! c_k^s (ks)_{q_{11}^{-1}}!=0, \]
a contradiction. Thus each $x_1^kx_2$ has infinite height. Thus,
$\left[x_1^kx_2\right]_c$ is a PBW-generator of infinite height.
If $r \leq n$ and $0 \leq n_1 \leq \dots \leq n_r, n_i \in \N$
such that $\sum_{j=1}^r n_j =n-r$, then $\left[x_1^{n_1}x_2\right]
\dots \left[x_1^{n_r}x_2\right]\in \toba^n(V)$ is an element of
$B_I'$. This collection is linearly independent, hence
    \[\dim \toba^n(V) \geq \sum^{n}_{r=1} \binom{ (n-r)+r-1}{ n-r }
    = \sum^{n-1}_{r=0} \binom{ n-1}{ r } =2^{n-1}. \]
Therefore, $\gkd (\toba(V))= \infty$. \epf

\section{The Weyl groupoid}\label{weyl}

\subsection{Groupoids}

\

There are several alternative definitions of a groupoid; let us
simply say that a groupoid is a category (whose collection of all
arrows is a set) where all the arrows are invertible. Let $\G$ be
a groupoid; it induces an equivalence relation $\approx$ on the
set of objects (or points) $\Pc$ by $x\approx y$ iff there exists
an arrow $g\in \G$ going from $x$ to $y$. If $x\in \Pc$, then
$\G(x) = $ all arrows going from $x$ to itself, is a group. A
groupoid is essentially determined by
\begin{itemize}
    \item the equivalence relation $\approx$, and
    \item the family of groups $\G(x)$, where $x$ runs in a set of
    representants of the partition associated to $\approx$.
\end{itemize}

A relevant example of a groupoid (for the purposes of this paper)
is the \emph{transformation groupoid}: if $G$ is a group acting on
a set $X$, then $\G = G\times X$ is a groupoid with operation
$(g,x)(h,y) = (gh, y)$ if $x=h(y)$, but undefined otherwise. Thus
the set of points in $\G$ is $X$ and an arrow $(g,x)$ goes from
$x$ to $g(x)$:
$$
\xymatrix{x \ar@/^/[0,2]^{(g,x)} & &   g(x) }
$$
In this example, $\G(x)$ is just the isotropy group of $x$. Thus,
if $G$ acts freely on $X$ (that is, all the isotropy groups are
trivially) then $\G$ is just the equivalence relation whose
classes are the orbits of the action. This is the case if
\begin{equation}\label{exa:transfgpd}
G= GL(\theta, \Z),\quad X = \text{ set of all ordered bases of }
\Z^\theta,
\end{equation}
with the natural action.

\subsection{The $i$-th reflection}

\

For $i,j \in \unon, i \neq j$, we consider
\begin{align*}
M_{i,j}&:=\left\{ (ad_cx_i)^m(x_j): m \in \N \right\}, \\
m_{ij} &:= \min \left\{ m \in \N:
(m+1)_{q_{ii}}(1-q_{ii}^mq_{ij}q_{ji})=0 \right\}.
\end{align*}
By Lemma \ref{conditions},  $m_{i,j} < \infty$ if and only if
$M_{i,j}$ is finite. In this case,
    \[(ad_cx_i)^{m_{ij}}(x_j)\neq 0, \ \ (ad_cx_i)^{m_{ij}+1}(x_j)=0. \]

Let $i \in \unon$. Set $m_{ii} = -2$. If any set $M_{i,j}$ is
finite, for all $ j \in \unon$, $i\neq j$, then we define a linear
map $s_i:\zt \to \zt$, by $s_i(e_j) =  e_j+m_{ij}e_i$,  $j\in
\unon$. Note that $s_i^2=\id$.

\bigbreak We recall that there are operators $y_i^L, y_i^R: \bB(V)
\rightarrow \bB(V)$, $i=1,\ldots ,n$ that play the role of left
and right invariant derivations. There is next a Hopf algebra
$H_i:= \kk \left[y_i^R\right] \# \kk \left[ e_i, e_i^{-1}
\right]$; $\toba(V)$ is an $H_i$-module algebra. We explicitly
record the following equality in $\cA_i := \left( \toba(V)^{op} \#
H_i^{cop} \right)^{op}$:
\begin{equation}\label{nose}
    \left( \rho \# 1\right)\pop \left( 1 \# y_i^R \right)
    = \left( 1 \# y_i^R \right)\pop (e_i^{-1} \td \rho) \# 1 + y_i^R(\rho) \# 1, \quad \rho \in \bB(V).
\end{equation}

In the setting above, the following Lemma is crucial for the proof
of Theorem \ref{transfnichols}. See \cite{Ang} for a complete
proof, slightly different from the argument sketched in \cite{H}.

\begin{lema}\label{keryil}
 $\bB(V) \cong \ker (y_i^L) \ot \kk \left[
x_i \right]$ as graded vector spaces. Moreover, $\ker (y_i^L)$ is
generated as algebra by $\cup _{j \neq i} M_{i,j}$. \qed
\end{lema}

The next result is the basic ingredient of the Weyl groupoid. We
discuss some details of the proof that are implicit in \cite{H}.

\begin{theorem}\label{transfnichols} \cite[Prop. 1]{H}.
Let $i \in \unon$ such that $M_{i,j}$ is finite, for all $ j \in
\unon$, $i\neq j$. Let $V_i$ be the vector subspace of $\cA_i$
generated by $\left\{ (ad_cx_i)^{m_{ij}}(x_j): j \neq i \right\}
\cup \left\{y_i^R\right\}$. The subalgebra $\bB_i$ of $\cA_i$
spanned by $V_i$ is isomorphic to $\bB(V_i)$, and
\[\Delta ^{+} (\bB_i) = \left\{ s_i \left( \Delta^{+} \left(
\bB(V) \right)\right) \setminus \left\{-e_i\right\} \right\} \cup
\left\{e_i\right\}. \]
\end{theorem}

\pf We just comment the last statement. The algebra $H_i$ is
$\zt$-graded, with $\deg y_i^R= -\e_i, \deg e_i^{\pm 1}=0$. Hence,
the algebra $\cA_i$ is $\zt$-graded, because $\bB(V)$ and $H_i$
are graded, and \eqref{nose} holds.

Hence, consider the abstract basis $\left\{u_j\right\}_{j \in
\unon}$ of $V_i$, with the grading $\deg u_j=\e_j$, $\bB(V_i)$ is
$\zt$-graded. Consider also the algebra homomorphism $\Omega:
\bB(V_i) \rightarrow \bB_i$ given by     \[ \Omega(u_j):= \left\{
\begin{array}{lc} (ad_c x_i)^{m_{ij}}
    (x_j) & j \neq i \\ y_i^R & j=i. \end{array} \right. \]
By the first part of the Theorem, proved in \cite{H}, $\Omega$ is
an isomorphism. Note:
\begin{itemize}
    \item $\deg \Omega (u_j)= \deg \left( \left( ad_c x_i \right)
    ^{m_{ij}} (x_j) \right) = \e_j+m_{ij}\e_i= s_i(\deg \ub_j)$,
    if $j \neq i$;
    \item $\deg \Omega (u_i)= \deg \left( y_i^R \right) = -\e_i= s_i(\deg
    \ub_i)$.
\end{itemize}
As $\Omega$ is an algebra homomorphism,  $\deg (\Omega(\ub))= s_i
(\deg(\ub))$, for all $\ub \in \bB(V_i)$. As $s_i^2=\id$, $s_i(
\deg (\Omega(\ub))) = \deg(\ub)$, for all $\ub \in \bB(V_i)$, and
$\hlb_{\bB(V_i)}= s_i (\hlb_{\bB_i})$.

\emph{Suppose first that $\ord x_i=h_i< \infty$.} Then
\[ \hlb_{\bB_i} = \hlb_{\ker y_i^L} \hlb_{\kk \left[y_i^R
\right]}= \hlb_{\ker y_i^L} \bq_{h_i}(X_i^{-1}) =
\frac{\hlb_{\bB(V)}}{\bq_{h_i}(X_i)}
 \bq_{h_i}(X_i^{-1}) , \]
the first equality because of $\Delta(\bB(V)) \subseteq
\N_0^{\theta}$, the second since $\ord x_i=\ord y_i^R$, and the
last by Proposition \ref{keryil}. As $s_i$ is an algebra
homomorphism, we have
\[ \hlb_{\bB(V_i)}= s_i (\hlb_{\bB_i}) = s_i \left(\hlb_{\bB(V)} \right) \frac{\bq_{h_i}(X_i)}{\bq_{h_i}(X_i^{-1})}. \]
But \begin{align*}
 s_i \left( \hlb_{\bB(V)} \right) &=
\prod_{\alpha \in \Delta^+ (\bB(V))} s_i \left(
\bq_{h_{\alpha}}(X^{\alpha}) \right) \\ &=  \left( \prod_{\alpha
\in \Delta^+ (\bB(V)) \setminus \{\e_i\}}
\bq_{h_{\alpha}}(X^{s_i(\alpha)}) \right) \bq_{h_i}(X_i^{-1});
\end{align*}
thus \[ \hlb_{\bB(V_i)} = \left( \prod_{\alpha \in \Delta^+
(\bB(V)) \setminus \{\e_i\}} \bq_{h_{\alpha}}(X^{s_i(\alpha)})
\right) \bq_{h_i}(X_i). \] By Lemma \ref{series}, $\Delta^{+}
(\bB_i) = \left\{ s_i \left( \Delta^{+} \left( \bB(V)
\right)\right) \setminus \left\{-e_i\right\} \right\} \cup
\left\{e_i\right\}$.

\emph{Suppose now that $\ord x_i=h_i = \infty$.} We have to
manipulate somehow the Hilbert series, because $\cA_i$ is not
locally finite. For this, we introduce an extra variable $X_0$,
corresponding to an  extra generator $\wte_0$ of $\Z^{\theta +1}$
(whose canonical  basis is denoted $\wte_0, \wte_1, \ldots,
\wte_{\theta}$), and consider $\Lambda = \frac{1}{2} \Z^{\theta
+1}$. We then define a $\Lambda$-grading in
$\bB(V)$, by \[ \widetilde{\deg}(x_j) = \begin{cases} \wte_j  &j\neq i,\\
    \frac12(\wte_i-\wte_0), \quad &j = i.\end{cases} \]

Let $\widetilde{s_i}:\Lambda \to \Lambda$ given by
 \[ \widetilde{s_i}(\wte_j) = \begin{cases} \wte_j + \frac{m_{ij}}{2}(\wte_i-\wte_0) &j\neq i,0,\\
    \wte_0, \quad &j = i, \\
    \wte_i, \quad &j = 0.\end{cases} \]
Consider the homomorphism $\Xi: \Lambda \rightarrow \zt$, given by
\[ \Xi(\widetilde{\e}_j) = \begin{cases} e_j &j\neq 0,
 \\ -e_i, \quad & j = 0.\end{cases} \]
Hence $\widetilde{\deg} (x_j)= \Xi (\deg x_j)$, for each $j \in
\unon$. Note that \begin{itemize}
    \item $\Xi (\widetilde{s_i}(\wte_j))= \Xi(\wte_j+ \frac{m_{ij}}{2}
    (\wte_i-\wte_0))= e_j+m_{ij}e_i=s_i(e_j)=s_i(\Xi(\wte_j))$, if     $j \neq i$,
    \item $\Xi (\widetilde{s_i}(\wte_i))= \Xi(\wte_0)= -e_i=s_i(e_i)=s_i(\Xi(\wte_i))$,
    \item $\Xi (\widetilde{s_i}(\wte_0))= \Xi(\wte_i)=
    e_i=s_i(-e_i)=s_i(\Xi(\wte_0))$;
\end{itemize}
thus $\Xi(\widetilde{s}_i(\alpha))= s_i (\Xi(\alpha))$, for all
$\alpha \in \Lambda$. With respect to grading, we can repeat the
previous argument, definining $\widetilde{\Delta}^+(\bB(V))
\subseteq \Lambda$ in analogous way to $\Delta^+(\bB(V))$. We get
\[ \widetilde{\Delta}^+ (\bB(V_i)) = \left(
\widetilde{s}_i(\widetilde{\Delta}^+ (\bB(V)) \setminus \{ -\wte_i
\} \right) \cup \{ \wte_i\}; \] as
$\Xi(\widetilde{\Delta}^+(\bB(V)))= \Delta^+(\bB(V))$,
$\Xi(\widetilde{\Delta}^+(\bB(V_i)))= \Delta^+(\bB(V_i))$, we have
\begin{eqnarray*}
\Delta^+(\bB(V_i)) &=& \Xi \left( \left( \widetilde{s}_i
(\widetilde{\Delta}^+ (\bB(V)) \setminus \{ -\wte_i \} \right)
\cup \{ \wte_i\} \right)\\ &=& s_i \left( \Xi \left(
\widetilde{\Delta}^+(\bB(V))\setminus \{\wte_i \} \right) \cup \{
e_i \} \right)
\\ &=& \left(
s_i \left( \Delta^{+} \left( \bB(V) \right)\right) \setminus
\left\{-e_i\right\} \right) \cup \left\{e_i\right\}.
\end{eqnarray*}
The proof now follows from Lemma \ref{series}. \epf

By Theorem \ref{transfnichols}, the initial braided vector space
with matrix $(q_{kj})_{1 \leq k,j \leq \theta}$,  is transformed
into another braided vector space of diagonal type $V_i$, with
matrix $(\overline{q}_{kj})_{1 \leq k,j \leq \theta}$, where
$\overline{q}_{jk}=q_{ii}^{m_{ij}m_{ik}}q_{ik}^{m_{ij}}q_{ji}^{m_{ik}}q_{jk}$,
$j,k \in\unon$.

If $j \neq i$, then $\overline{m_{ij}} = \min \left\{ m \in \N :
(m+1)_{\overline{q}_{ii}}
    \left(\overline{q}_{ii}^m\overline{q}_{ij} \overline{q}_{ji}=0 \right) \right\} \overset{\!}= m_{ij}$.
Thus, the previous transformation is invertible.

\subsection{Definition of the Weyl groupoid}

\

Let $E = (\e_1, \dots, \e_\theta)$ be the canonical basis of
$\Z^\theta$. Let $(q_{ij})_{\otvz}\in
(\C^{\times})^{\theta\times\theta}$. We fix once and for all the
bilinear form $\chi: \Z^\theta \times \Z^\theta \to\C^{\times}$
given by
$$
\chi(\e_i, \e_j) = q_{ij}, \qquad \otvz.$$

Let $F= (\f_1, \dots, \f_\theta)$ be an arbitrary ordered basis of
$\Z^\theta$ and let $\qf_{ij} = \chi(\f_i,\f_j)$, $\otvz$, the
\emph{braiding matrix with respect to the basis $F$}. Fix $i\in
\unon$. If $\otrvz$, then we consider the set
$$
\{m\in \N_0: (m+1)_{\qf_{ii}} \, (\qf_{ii}^m\qf_{ij}\qf_{ji} - 1)
= 0 \}.
$$
This set might well be empty, for instance if $\qf_{ii} = 1 \neq
\qf_{ij}\qf_{ji}$ for all $j\neq i$. If this set is nonempty, then
its minimal element is denoted $\mf_{ij}$ (which of course depends
on the basis $F$). Set also $\mf_{ii} = 2$ . Let $\si\in
GL(\Z^{\theta})$ be the pseudo-reflection given by $\si(\f_j) :=
\f_j + m_{ij}\f_i, \quad j\in \unon. $

Let us compute the braiding matrix with respect to the matrix
$\si(F)$. Let $\ub_j := \si(\f_j)$ and $\qu_{rj} = \chi(\ub_r,
\ub_j)$. If we also set $\mf_{ii} := -2$, $\otv$ by convenience,
then
\begin{equation}\label{eqn:braiding-matrix}
\qu_{rs} = {\qf_{ii}}^{\;m_{ir}m_{is}}{\qf_{ri}}^{\;m_{is}}
{\qf_{is}}^{\;m_{ir}}{\qf_{rs}}, \qquad 1\le r, s \le \theta.
\end{equation}
In particular $\qu_{ii} = \qf_{ii}$ and $\qu_{jj} =
\left(\qf_{ii}^{\;m_{ij}}\qf_{ji}\qf_{ij}\right)^{\;m_{ij}}
{\qf_{jj}}$, $1\le j \le \theta$. Thus, even if $m_{ij}$ are
defined for the basis $F$ and for all $i\neq j$, they need not be
defined for the basis $\si(F)$. For example, if
$$
\begin{pmatrix}
  \qf_{ii} & \qf_{ij} \\ \qf_{ji} & \qf_{jj} \\
\end{pmatrix} = \begin{pmatrix}
  -1 & -\xi \\ \xi & \xi^{-2} \\
\end{pmatrix},
$$
where $\xi$ is a root of 1 of order $>4$, then $\begin{pmatrix}
  \qu_{ii} & \qu_{ij} \\ \qu_{ji} & \qu_{jj} \\
\end{pmatrix} = \begin{pmatrix}
  -1 & \xi^{-1} \\ -\xi^{-1} & 1 \\
\end{pmatrix}$ and
$m_{ji}$ is not defined with respect to the new basis $\si(F)$.
However, for an arbitrary $F$ and $i$ such that $m_{ij}$ for $F$
is defined, then $m_{ij}$ is defined for the new basis $\si(F)$
and coincides with the old one, so that
\begin{equation}\label{eqn:si-old}
s_{i, \si(F)} = \si.
\end{equation}

Notice that formula \eqref{eqn:braiding-matrix} and a variation
thereof appear in \cite{H2}.

\begin{definition}\label{defi:weyl-groupoid}
The Weyl groupoid $W(\chi)$ of the bilinear form $\chi$ is the
smallest subgroupoid of the transformation groupoid
\eqref{exa:transfgpd} with respect to the following properties:
\begin{itemize}
    \item $(\id, E)\in W(\chi)$,
    \item if $(\id, F)\in W(\chi)$ and $\si$ is defined, then $(\si, F)\in
    W(\chi)$.
\end{itemize}
\end{definition}

In other words, $W(\chi)$ is just a set of bases of $\Z^\theta$:
the canonical basis $E$, then all bases $s_{i, E}(E)$ provided
that $s_{i,E}$ is defined, then all bases $s_{j, s_i(E)}s_{i,
E}(E)$ provided that $s_{i,E}$ and $s_{j, s_i(E)}$ are defined,
and so on.

\bigbreak Here is an alternative description of the Weyl groupoid.
Consider the set of all pairs $(F, (\qf_{ij})_{\otvz})$ where $F$
is an ordered basis of $\Z^\theta$ and the $\qf_{ij}$'s are
non-zero scalars. Let us say that $$(F, (\qf_{ij})_{\otvz})\sim
(U, (\qu_{ij})_{\otvz})$$ if there exists and index $i$ such that
$m_{ij}$ exists for all $\otrvz$, $U = \si(F)$ and $\qu_{ij}$ is
obtained from $\qf_{ij}$ by \eqref{eqn:braiding-matrix}. By
\eqref{eqn:si-old} this is reflexive; consider the equivalence
relation $\approx$ generated by $\sim$. Then $W(\chi)$ is the
equivalence class of $(E, (q_{ij})_{\otvz})$ with respect to
$\approx$. Actually, if $(F, (\qf_{ij})_{\otvz})$ and $(U,
(\qu_{ij})_{\otvz})$ belong to the equivalence class, there will
be a unique $s\in GL(\Z^\theta)$, which is a product of suitable
$s_i$'s, such that $(s,F)\in W(\chi)$ and $s(F) = U$.

The equivalence class of $(F, (\qf_{ij})_{\otvz})$ is denoted
$W(F, (\qf_{ij})_{\otvz})$. Furthermore, if $\chi$ is a fixed
bilinear form as above, $F = (f_i)_{\otv}$ and $\qf_{ij} =
\chi(f_i,f_j)$, $\otvz$, then we denote $W(F, \chi) := W(F,
(\qf_{ij})_{\otvz})$; and $W(\chi) := W(E, \chi)$ where $E$ is the
canonical basis.

From this viewpoint, it is  natural to introduce the following
concept.

\begin{definition}\label{defi:weyl-equivalence}
We say that $(\qf_{ij})_{\otvz}$ and $(\qu_{ij})_{\otvz} \in
\kk^{\times}$ are \emph{Weyl equivalent} if there exist ordered
bases $F$ and $U$ such that
$$(F, (\qf_{ij})_{\otvz})\approx (U, (\qu_{ij})_{\otvz}).$$
Now recall that $(\qf_{ij})_{\otvz}$ and $(\qu_{ij})_{\otvz} \in
\kk^{\times}$ are \emph{twist equivalent} if $\qf_{ii} = \qu_{ii}$
and $\qf_{ij}\qf_{ji} = \qu_{ij}\qu_{ji}$ for all $\otvz$.

It turns out that it is natural to consider the equivalence
relation $\wh$ generated by twist- and Weyl-equivalence, see
\cite[Def. 3]{H2}, see also \cite[Def. 2]{H3}. We propose to call
$\wh$ the \emph{Weyl-Heckenberger equivalence}; note that this is
the ``Weyl equivalence" in \cite{H2}. We suggest this new
terminology because the Weyl groupoid is really an equivalence
relation.
\end{definition}

\bigbreak The Weyl groupoid is meant to generalize the set of
basis of a root system. For convenience we set $\Bg(\chi) = \{F:
(\id, F)\in W(\chi)\}$, the set of points of the groupoid
$W(\chi)$. Then the \emph{generalized root
system}\footnote{Actually this is a little misleading, since in
the case of braidings of symmetrizable Cartan type, one would get
just the real roots. } associated to $\chi$ is
\begin{equation}\label{eqn:root-system}
\Delta(\chi) = \bigcup_{F\in \Bg(\chi)} F.
\end{equation}

We record for later use the following evident fact.

\begin{obs}\label{obs:gpd-finite}
The following are equivalent:
\begin{enumerate}
    \item The groupoid $W(\chi)$ is finite.
    \item The set $\Bg(\chi)$ is finite.
    \item The generalized root system $\Delta(\chi)$ is finite.
    \qed
\end{enumerate}
\end{obs}

Let also $\map: W(\chi) \to GL(\theta, \Z)$, $\map(s,F) = s$ if
$(s, F)\in W(\chi)$. We denote by $\wo$ the subgroup generated by
the image of $\map$. Compare with \cite{serganova}.

\bigbreak Let us say that $\chi$ is \emph{standard} if for any
$F\in \Bg(\chi)$, the integers $m_{rj}$ are defined, for all $1\le
r,j\le \theta$, and the integers $m_{rj}$ for the bases $\si(F)$
coincide with those for $F$ for all $i,r,j$ (clearly it is enough
to assume this for the canonical basis $E$).

\begin{prop}\label{prop:standard}
Assume that $\chi$ is standard. Then $$\wo = \langle \sE:
\otv\rangle.$$ Furthermore $\wo$ acts freely and transitively on
$\Bg(\chi)$. \qed
\end{prop}
The first claim says that $\wo$ is a Coxeter group. The second
implies that $\wo$ and $\Bg(\chi)$ have the same cardinal.

\pf Let $F\in \Bg(\chi)$. Since $\chi$ is standard, for any
$\otvz$
\begin{equation}\label{eqn:si-standard}
s_{j, \si(F)} = \si s_{j, F}\si.
\end{equation}
Hence $\wo \subseteq \langle \sE: \otv\rangle$; the other
inclusion being clear, the first claim is established. Now, by the
very definition of the Weyl groupoid, there exists a unique $w\in
\wo$ such that $w(E) = F$. Thus, to prove the second claim we only
need to check that the action is well-defined; and for this, it is
enough to prove: if $w\in \wo$, then $w(E)\in \Bg(\chi)$. We
proceed by induction on the length of $w$, the case of length one
being obvious. Let $w' =  w\sE$, with length of $w' =$ length of
$w + 1$. Then $F = w(E)\in \Bg(\chi)$. the matrix of $\si$ in the
basis $E$ is $\| \si\|_{E} = \|\id\|_{F,E}\| \si\|_{F}
\|\id\|_{E,F}$ and since $\chi$ is standard, we conclude that $\si
= w\sE w^{-1}$. \footnote{Here, one uses that the matrix
$\|\id\|_{F,E}$ when seen as transformation of $\Z^\theta$, sends
$\e_i$ to $\f_i$ for all $i$.} That is, $w' = \si w$ and $w'(E) =
\si(F)\in \Bg(\chi)$. \epf

\begin{obs}\label{obs:gpd-finite-standard}
Assume that $\chi$ is standard. Then the following are equivalent:
\begin{enumerate}
    \item The groupoid $W(\chi)$ is finite.
    \item The set $\Bg(\chi)$ is finite.
    \item The generalized root system $\Delta(\chi)$ is finite.
    \item The group $\wo$ is finite.
\end{enumerate}
If this holds, then $\wo$ is a finite Coxeter group; and thus
belongs to the well-known classification list in \cite{B}.
\end{obs}

\subsection{Nichols algebras of Cartan type}

\

\begin{definition}
A braided vector space $(V,c)$ is of \emph{Cartan type} if it is
of diagonal type with matrix $(q_{ij})_{1 \leq i,j \leq \theta}$
and for any $i,j \in \unon$,  $q_{ii} \neq 1$, and there exists
$a_{ij} \in \Z$ such that
    \[q_{ij}q_{ji}= q_{ii}^{a_{ij}}. \]
\end{definition}
\noindent The integers $a_{ij}$ are uniquely determined by
requiring that
 $a_{ii}=2$, $0 \leq -a_{ij} < \ord(q_{ii})$, $\otrvz$.
Thus  $(a_{ij})_{1 \leq i,j \leq \theta}$ is a generalized Cartan
matrix \cite{K}.

\medbreak If $(V,c)$ is a braided vector space of Cartan type with
generalized Cartan matrix $(a_{ij})_{1 \leq i,j \leq \theta}$,
then for any $i,j \in \unon, j \neq i$, $m_{ij}=-a_{ij}$. It is
easy to see that a braiding of Cartan type is standard, see the
first part of the proof of \cite[Th. 1]{H}. Hence we have from
Remark \ref{obs:gpd-finite-standard}:

\begin{lema}\label{lema:cartan-finite}
Assume that $\chi$ is of Cartan type with symmetrizable Cartan
matrix $C$. Then the following are equivalent:
\begin{enumerate}
    \item The generalized root system $\Delta(\chi)$ is finite.
    \item The Cartan matrix $C$ is of finite type. \qed
\end{enumerate}
\end{lema}

We are now ready to sketch a proof of the main theorem in
\cite{H}.

\begin{theorem}\label{heckenberger}
Let $V$ be a braided vector space of Cartan type with Cartan
matrix $C$. Then, the following are equivalent.
\begin{enumerate}
    \item The set $\Delta(\bB(V))$ is finite.
    \item The Cartan matrix $C$ is of finite type.
\end{enumerate}
\end{theorem}
\pf (1) $\Rightarrow$ (2). As $\dc \subseteq \Delta(\bB(V))$,
$\dc$ is finite. If $C$ is symmetrizable, we apply Lemma
\ref{lema:cartan-finite}. If $C$ is not symmetrizable, one reduces
as in \cite{AS-adv} to the smallest possible cases, see \cite{H}.
See \cite{H} for the proof of (2) $\Rightarrow$ (1). \epf

We can now prove Lemma \ref{rosso}: (b) $\implies$ (a) was already
discussed in \cite{AS-crelle}. (a) $\implies$ (b): it follows from
Lemma \ref{rossito} that $(V,c)$ is of Cartan type; hence it is of
finite Cartan type, by Theorem \ref{heckenberger}. To prove that
$(V,c)$ is of DJ-type-- see \cite[p. 84]{AS-crelle}-- it is enough
to assume that $C$ is irreducible; then the result follows by
inspection.

\bigbreak We readily get the following Corollary, as in \cite[Th.
2.9]{AS-crelle}-- that really follows from results of Lusztig and
Rosso. Let $V$ be a braided vector space of diagonal type with
matrix $q_{ij} \in \kk^{\times}$. Let $m_{ij}\geq 0$ be as in
Lemma \ref{rossito}.

\begin{cor}\label{quantumserre}
If $\toba(V)$ is a domain and its Gelfand-Kirillov dimension is
finite, then $\toba(V) \simeq \kk <x_1, \dots, x_\theta:
(ad_cx_i)^{m_{ij} + 1}(x_j)=0, \quad i\neq j >$. \qed
\end{cor}

Notice that the hypothesis ``$\toba(V)$ is a domain" is equivalent
to ``$q_{ii} = 1$ or it is not a root of 1, for all $i$",
\emph{cf.} \cite{Kh}.

\appendix

\section{Generic data and the definition of $U(\mathcal
D)$}

In this Appendix, we briefly recall the main definitions and
results from \cite{AS-crelle} needed for Theorem
\ref{fingrowth-lifting}. Everything below belongs to
\cite{AS-crelle}; see \emph{loc. cit.} for more details. Below, we
shall refer to the following terminology.

\begin{itemize} \item $\Gamma$ is a free abelian group of finite rank
$s$.

\item $(a_{ij})\in \Z^{\theta\times
\theta}$ is a Cartan matrix  of finite type \cite{K}; we denote by
$(d_{1}, \dots, d_{\theta})$ a diagonal matrix of positive
integers such that $d_{i}a_{ij} = d_{j} a_{ji}$, which is minimal
with this property.

\item $\mathcal X$ is the set of connected components of the
Dynkin diagram corresponding to the Cartan matrix $(a_{ij})$; if
$i,j \in \{1, \dots, \theta\}$, then $i\sim j$ means that they
belong to the same connected component.

\item $(q_{I})_{I\in \mathcal X}$ is a  family of elements
in $\ku$ which are not roots of 1.

\item $g_{1}, \dots, g_{\theta}$ are elements in $\Gamma$,
$\chi_{1}, \dots, \chi_{\theta}$ are characters in
$\widehat{\Gamma}$, and all these satisfy
\begin{equation}
\label{cartantype} \langle \chi_{i}, g_{i}\rangle = q_{I}^{d_i},
\quad \langle \chi_{j}, g_{i}\rangle \langle \chi_{i},
g_{j}\rangle = q_{I}^{d_{i}a_{ij}}, \quad \text{for all }  1 \le
i, j \le \theta, \quad i\in I. \notag
\end{equation}
\end{itemize}

We say that two vertices $i$ and $j$ {\em are linkable} (or that
$i$ {\em is linkable to} $j$) if $i\not\sim j$, $g_{i}g_{j} \neq
1$ and $\chi_{i}\chi_{j} = \varepsilon$.

\begin{definition}\label{link-dat}
A {\em linking datum}  for
$$\Gamma, \quad (a_{ij}), \quad(q_{I})_{I\in \mathcal X},
\quad g_{1}, \dots, g_{\theta}\quad \text{and } \chi_{1}, \dots,
\chi_{\theta}$$ is a collection $(\lambda_{ij})_{1 \le i < j \le
\theta, i \nsim j}$ of elements in $\{0, 1\}$ such that
$\lambda_{ij}$ is arbitrary if $i$ and $j$ are linkable but 0
otherwise. Given a linking datum, we say that two vertices $i$ and
$j$ {\em are linked} if $\lambda_{ij}\neq 0$. The collection
$$\mathcal D = \mathcal D((a_{ij}),(q_{I}), (g_i), (\chi_i),
(\lambda_{ij})),$$ where $(\lambda_{ij})$ is a linking datum, will
be called a {\it generic datum of finite Cartan type} for
$\Gamma$.
\end{definition}

In the next Definition, $\ad$ is the ``braided" adjoint
representation, see \cite{AS-crelle}.

\begin{definition}\label{ud}
Let $\mathcal D = \mathcal D((a_{ij}),(q_{I}), (g_i), (\chi_i),
(\lambda_{ij}))$ be a generic datum of finite Cartan type for
$\Gamma$. Let $U({\mathcal D})$ be the algebra presented by
generators $a_{1}, \dots, a_{\theta}$, $y_{1}^{\pm 1}, \dots,
y_{s}^{\pm 1}$ and relations
\begin{align}\label{relations}
 y_{m}^{\pm 1}y_{h}^{\pm 1} &= y_{h}^{\pm 1}y_{m}^{\pm 1}, \quad y_{m}^{\pm 1}
y_{m}^{\mp 1} = 1, \qquad 1 \le m,h  \le s,\\\label{relations1}
y_{h}a_{j} &= \chi_{j}(y_{h})a_{j}y_{h}, \qquad 1 \le h \le s,\, 1 \le j \le \theta,\\
\label{relations2}
(\ad a_{i})^{1 - a_{ij}}a_{j} &= 0, \qquad 1 \le i \neq j \le \theta, \quad i \sim j, \\
\label{relations3} a_{i}a_{j} - \chi_{j}(g_{i})a_{j}a_{i} &=
\lambda_{ij}(1 - g_{i}g_{j}), \qquad 1 \le i < j \le \theta, \quad
i \not\sim j.
\end{align}
\end{definition}

The relevant properties of $U({\mathcal D})$ are stated in the
following result.

\begin{theorem}\label{construction} \cite[Th. 4.3]{AS-crelle}.
The algebra $U({\mathcal D})$ is a pointed Hopf algebra with
structure determined  by
\begin{equation}\label{sk-gl}
\Delta y_{h} = y_{h}\otimes y_{h}, \qquad \Delta a_{i} =
a_{i}\otimes 1 + g_{i} \otimes a_{i}, \qquad 1 \le h \le s,\, 1
\le i \le \theta.
\end{equation}

Furthermore,  $U({\mathcal D})$ has a PBW-basis given by monomials
in the root vectors, that are defined by an iterative procedure.
The coradical filtration of $U({\mathcal D})$ is given by the
ascending filtration in powers of those root vectors. The
associated graded Hopf algebra $\gr U({\mathcal D})$ is isomorphic
to $\toba(V)\# \ku \Gamma$; $U({\mathcal D})$ is a domain with
finite Gelfand-Kirillov dimension.
\end{theorem}

\subsection*{Acknowledgements} We thank Istv\'an Heckenberger for
some useful remarks on a first version of this paper. We also
thank the referee for his/her interesting remarks.

\end{document}